\documentclass[11pt,a4paper]{article}%
\usepackage{amsmath}
\usepackage{amsfonts}
\usepackage{amssymb}
\usepackage{graphicx}%
\begin{document}

\title{APPLICATIONS OF THE JACOBI GROUP TO QUANTUM MECHANICS}
\author{S. BERCEANU\,$^{1}$, A. GHEORGHE\,$^{1}$\\$^{1}$ National Institute for Nuclear Physics and Engineering,\\P.O.Box MG-6, RO-077125 Bucharest-Magurele, Romania,\\E-mail: Berceanu@theory.nipne.ro}
\maketitle

\begin{abstract}
Infinitesimal holomorphic realizations for the Schr\"{o}dinger-Weil
representation and the discrete series representations of the Jacobi group are
constructed. Explicit expressions of the basic differential operators are
obtained. The squeezed states for the unitary irreducible representation of
the Jacobi group are introduced. Matrix elements of the squeezed operators,
expectation values of polynomial operators in infinitesimal generators of the
Jacobi group, the squeezing region and a description of Mandel's parameter are presented.

\end{abstract}


\section{Introduction}

The Jacobi group $G_{n}^{J}$ is the semidirect product of the symplectic group
$Sp(n,\mathbb{R})$ with an appropriate Heisenberg group \cite{ez,bs,yang}. In
\cite{jac1,sbj,sbc} we have considered the Perelomov coherent states and the
squeezed states for the Jacobi group. The representation theory of the Jacobi
group has been constructed in \cite{bs,yang,bebo,tak,taka2} with relevant
topics: Schr\"{o}dinger-Weil and metaplectic subrepresentations,
classification and realizations of irreducible unitary representations over
local fields, holomorphic Jacobi forms and automorphic representations,
symplectic orbits, Whittaker models, Hecke and generalized Kac-Moody algebras,
$L$-functions and modular forms, spherical functions, and the ring of
invariant differential operators.

In section 2.1 we review briefly some basic facts about the Jacobi group
$G^{J}=$ $G_{1}^{J}$. We realize the infinitesimal generators of $G^{J}$ in
terms of boson operators and standard infinitesimal generators of the
symplectic group. In \ section 2.2 we construct infinitesimal holomorphic
realizations for the Schr\"{o}dinger-Weil representation and the discrete
series representations. In section 3 we introduce appropriate squeezed states
for these representations.

\section{Unitary representations of the Jacobi group}

\subsection{The fundamental principle}

The real Jacobi group $G^{J}$ is a subgroup of the symplectic group
$Sp(2,\mathbb{R)}$ consisting of $4\times4$ real matrices $g=\left(  \left(
\lambda,\mu,\kappa\right)  ,M\right)  $ of the form \
\begin{equation}
g=\left(
\begin{array}
[c]{cccc}%
a & 0 & b & a\mu-b\lambda\\
\lambda & 1 & \mu & \kappa\\
c & 0 & d & c\mu-d\lambda\\
0 & 0 & 0 & 1
\end{array}
\right)  ,\ M=\left(
\begin{array}
[c]{cc}%
a & b\\
c & d
\end{array}
\right)  ,\ ad-bc=1,\ \label{defj}%
\end{equation}
where $\left(  \lambda,\mu,\kappa\right)  \in H(\mathbb{R)}$ and$\ M\in
SL(2,\mathbb{R)}$ \cite{bs}. \ Here $H(\mathbb{R)}$ is the three-dimensional
real Heisenberg group and $SL(2,\mathbb{R)=}Sp(1,\mathbb{R)}$ is the special
linear group. Then $G^{J}$ is the semidirect product of $H(\mathbb{R)}$ with
$SL(2,\mathbb{R)}$. Let $\mathfrak{\ }\left\langle X_{1},...,X_{n}%
\right\rangle _{\mathbb{F}}$ denote the Lie algebra over $\mathbb{F}$ \ with
the basis elements $X_{1},...,X_{n}$. We denote by $\mathbb{R}$, $\mathbb{C}$,
$\mathbb{Z}$, and $\mathbb{N}$ the field of real numbers, the field of complex
numbers, the ring of integers, and the set of non-negative integers,
respectively. The Lie algebras of $G^{J}$, $H(\mathbb{R)}$, and
$SL(2,\mathbb{R)}$ are denoted by $\mathfrak{g}^{J}=\left\langle
\mathsf{P,Q,R,F,G,H}\right\rangle _{\mathbb{R}}$, ${\mathfrak{h}}=\left\langle
\mathsf{P,Q,R}\right\rangle _{\mathbb{R}}$, ${\mathfrak{sl}}(2,\mathbb{R)}%
=\left\langle \mathsf{F,G,H}\right\rangle _{\mathbb{R}}$, respectively.
$\mathsf{P}$\textsf{, }$\mathsf{Q}$\textsf{, }$\mathsf{R}$\textsf{,
}$\mathsf{F}$\textsf{, }$\mathsf{G}$\textsf{, }$\mathsf{H}$ are $4\times4$
matrices of coefficients $\ \mathsf{F}_{ij}=\delta_{i1}\delta_{j3}%
$,$\ \mathsf{G}_{ij}=\delta_{i3}\delta_{j1}$, $\ \mathsf{H}_{ij}=\delta
_{i1}\delta_{j1}-\delta_{i3}\delta_{j3}$, $\mathsf{\ P}_{ij}=\delta_{i2}%
\delta_{j1}-\delta_{i3}\delta_{j4}$, $\ \mathsf{Q}_{ij}=\delta_{i1}\delta
_{j4}+\delta_{i2}\delta_{j3}$, $\ \mathsf{R}_{ij}=\delta_{i2}\delta_{j4}$,
where $i,j=1,2,3,4$. We get the commutators \cite{bs}
\begin{align}
\left[  \mathsf{P,Q}\right]   &  =\mathsf{2R},\mathsf{\ }\left[
\mathsf{F,G}\right]  =\mathsf{H},\mathsf{\ }\left[  \mathsf{H,F}\right]
=\mathsf{2F},\mathsf{\ }\left[  \mathsf{G,H}\right]  =\mathsf{2G}%
,\label{commj}\\
\lbrack\mathsf{P,F}]  &  =\mathsf{Q},\mathsf{\ }[\mathsf{Q,G}]=\mathsf{\ P}%
,\mathsf{\ }\left[  \mathsf{P,H}\right]  =\mathsf{\ P},\mathsf{\ }\left[
\mathsf{H,Q}\right]  =\ \mathsf{\ Q},\nonumber
\end{align}
all other are zero. The center of $G^{J}$, consisting of all $\left(
0,0,\kappa\right)  \in H(\mathbb{R)}$, is thus isomorphic to $\mathbb{R}$.
Every non-trivial central character $\psi$ of index $m\in\mathbb{R}$ can be
obtained as $\psi(\left(  0,0,\kappa\right)  )=\exp(2\pi\mathrm{i}m\kappa)$,
where $\left(  0,0,\kappa\right)  \in H(\mathbb{R)}$.

Let $\pi$ an unitary irreducible representation of $G^{J}$ of nonzero index
$m$ on a complex separable Hilbert space $\mathcal{H}$. Let $\hat{\pi}$ be the
derived representation of $\pi$ and let $\mathcal{D}$ be the space of smooth
vectors. Denote $X=\hat{\pi}(\mathsf{X})$ for any $\mathsf{X}\in
\mathfrak{g}^{J}$. Then $\hat{\pi}\left(  \mathsf{R}\right)  =\mathrm{i}\mu
I$, where $\mu=2\pi m$ and $I$ is the identity operator. We now introduce the
following operators in the complexification $\hat{\pi}(\mathfrak{g}%
_{\mathbb{C}}^{J}{\mathfrak{)}}$ of $\hat{\pi}(\mathfrak{g}^{J}{\mathfrak{)}}%
$:
\begin{align}
a  &  =\frac{1}{2\sqrt{\left\vert \mu\right\vert }}\left(
P-\mathbf{\mathrm{i}}\sigma Q\right)  ,\ a\mathbf{^{\dagger}}=-\frac{1}%
{2\sqrt{\left\vert \mu\right\vert }}\left(  P+\mathbf{\mathrm{i}}\sigma
Q\right)  ,\ \label{aK}\\
K_{\pm}  &  =\mp\frac{1}{2}H-\frac{\mathbf{\mathrm{i}}\sigma}{2}\left(
F+G\right)  ,\ K_{0}=\frac{\mathbf{\mathrm{i}}\sigma}{2}\left(  G-F\right)
\mathbf{,}\nonumber
\end{align}
where $\sigma=\mu/\left\vert \mu\right\vert $. Using (\ref{commj}) and
$X^{\dagger}=-X$ for $X\in\hat{\pi}(\mathfrak{g}^{J}{\mathfrak{)}}$, we obtain

\textbf{Proposition 1.} $\hat{\pi}(\mathfrak{g}_{\mathbb{C}}^{J}{\mathfrak{)}%
}=\left\langle I,a,a^{\dagger},K_{+},K_{-},K_{0}\right\rangle _{\mathbb{C}}$
\textit{with }$\left(  a^{\dagger}\right)  ^{\dagger}=a$,$\ $

\noindent$K_{\pm}^{\dagger}=K_{\mp}$, $\ \left[  a,{K}_{+}\right]
=a^{\dagger}$, $\left[  {K}_{-},a\right]  =0$, $2\left[  a,{K}_{0}\right]
=a$, \textit{and} \textit{ }%
\begin{equation}
\lbrack a,a^{\dagger}]=I\mathbf{,\ [}K_{0},K_{\pm}]=\pm K_{\pm}\ ,\ [K_{-}%
,K_{+}]=2K_{0}. \label{boson}%
\end{equation}
Then $\hat{\pi}({\mathfrak{sl}}(2,\mathbb{R))}=\left\langle K_{0},K_{1}%
,K_{2}\right\rangle _{\mathbb{R}}$, where $K_{\pm}=K_{1}\pm iK_{2}$. Consider
the following operators in the universal enveloping algebra $\mathcal{U(}%
\hat{\pi}(\mathfrak{g}_{\mathbb{C}}^{J}{\mathfrak{)}})$ \cite{bs}:%
\begin{equation}
W_{-}=K_{-}-\frac{1}{2}a^{2},~W_{+}=K_{+}-\frac{1}{2}(a^{\dagger})^{2}%
,\ W_{0}=K_{0}-\frac{1}{2}a^{\dagger}a-\frac{1}{4}I. \label{wald}%
\end{equation}
We have $\left[  W_{0},W_{\pm}\right]  =\pm W_{\pm}$,$\ \ \left[  W_{-}%
,W_{+}\right]  =2W_{0}$, $W_{\pm}^{\dagger}=W_{\mp}$ and $[a,W_{\sigma}]=0$,
where $\sigma=0,+,-$ . Let ${\mathfrak{w}}=\left\langle W_{0},W_{1}%
,W_{2}\right\rangle _{\mathbb{R}}$, where $W_{\pm}=W_{1}\pm\mathrm{i}W_{2}$.
The Casimir operator of ${\mathfrak{w}}\simeq{\mathfrak{sl}}(2,\mathbb{R)}$ is
defined by $C=W_{0}^{2}-W_{1}^{2}-$ $W_{2}^{2}$ . The metaplectic group
$Mp(2,\mathbb{R)}$ is the non-split two-fold cover of $SL(2,\mathbb{R)}$.

Using the Lie algebra ${\mathfrak{w}}$, the representation theory of $G^{J}$
may be fully reduced to Waldspurger's representation theory of
$Mp(2,\mathbb{R)}$ \cite{bs,wald}. The Stone-von Neumann theorem and the
method of Mackey for semidirect products lead us to the \emph{fundamental
principle\ }in the representation theory of the Jacobi group \cite{bs}:

Any representation $\pi$ of $G^{J}$ with index $m\neq0$ is obtained in a
unique way as $\pi=\pi_{\mathrm{SW}}^{m}\otimes\tilde{\pi}$, where the
\emph{Schr\"{o}dinger-Weil representation} $\pi_{\text{\textrm{SW}}}^{m}$ is a
certain projective representation of $G^{J}$ and $\tilde{\pi}$ is a
representation of the metaplectic group $Mp(2,\mathbb{R)}$ (considered as a
projective representation of $SL(2,\mathbb{R)}$). The representations $\pi$
and $\tilde{\pi}$ are simultaneously unitary, and irreducible.

The irreducible unitary representation $\pi$ of $G^{J}$ of index $m\neq0$ are
infinitesimally equivalent to the \emph{principal series} representations
$\hat{\pi}_{ms\nu}$ for $s\in\mathrm{i}\mathbb{R}\cup(-1/2,1/2)$, $\nu=\pm
1/2$, with $C=(s^{2}-1)/4$, or to the \emph{positive and negative discrete
series} representations $\hat{\pi}_{mk}^{\pm}$ for $k\in\mathbb{Z}$, $k\geq1$,
with $C=(k-1/2)(k-5/2)/4$ \cite{bs}. Here $2K_{0}$ has the integral dominants
weights $k$ for $\hat{\pi}_{mk}^{+}$ and $1-k$ for $\hat{\pi}_{mk}^{-}$.

Standard models of the preceding representations are presented in \cite{bs}.
Let $L_{\mathrm{hol}}^{2}(\mathcal{M},\nu)$ denote the complex Hilbert space
of all $\mathbb{C}$-valued holomorphic functions on the complex manifold
$\mathcal{M}$ which are square-integrable with respect to the measure $\nu$.
Consider the upper half-plane $\mathbb{H}$ of all $\tau\in\mathbb{C}$ with
$\Im\tau>0$. Suppose $k>3/2$ and $m\neq0$. Then there is an irreducible
unitary representation $\pi_{mk}$ on $L_{\mathrm{hol}}^{2}(\mathbb{C\times
H},\mu_{mk})$ with the Petersson measure $\mathrm{d}\mu_{mk}=\exp\left(  -4\pi
my^{2}/v\right)  v^{k-3}\mathrm{d}^{2}z\ \mathrm{d}^{2}\tau$, where $y=\Im z$
and $v=\Im\tau$ \cite{bs,bebo,taka2}. For any $f\in L_{\mathrm{hol}}%
^{2}(\mathbb{C\times H},\mu_{mk})$ and $g=\left(  \left(  \lambda,\mu
,\kappa\right)  ,M\right)  \in G^{J}$, we have
\begin{equation}
\pi_{mk}(g^{-1})f(z,\tau)=(c\tau+d)^{-k}\exp\left[  2\pi\mathrm{i}m\left(
\kappa+\theta\right)  \right]  f\left(  z_{g},\tau_{g}\right)  ,
\label{petersson}%
\end{equation}
where $(z,\tau)\in\mathbb{C\times H}$, $z_{g}=(c\tau+d)^{-1}(z+\lambda\tau
+\mu)$,$\ \tau_{g}=(c\tau+d)^{-1}(a\tau+b)$, and \noindent$\theta=\lambda
z+(\lambda z_{g}-cz_{g}^{2})(c\tau+d)$. The Jacobi forms \cite{ez} are
associated with $\pi_{mk}$ provided the index $m$ and the weight $k$ are
positive integers. $\pi_{mk}$ can also be used to produce bases for signal
processing, continuous windowed Fourier and wavelet transforms \cite{karen}.
The group $G^{J}$ is unimodular. The representation $\pi_{mk}$ is
square-integrable modulo the center of $G^{J}$. Then $\pi_{mk}$ is a Perelomov
coherent state representation based on $\mathbb{C\times H}$ \cite{karen}. Let
$\mathbb{D}$ be the open disk of the points $w\in\mathbb{C}$ with $\left\vert
w\right\vert <1$. The manifolds $\mathbb{C\times H}$ and $\mathbb{C\times D}$
can be biholomorphically identified by a partial Cayley transform \cite{bs}.
Then there is an irreducible unitary representation $\rho_{mk}$ on
$L_{\mathrm{hol}}^{2}(\mathbb{C\times D},\nu_{mk})$ which is unitarily
equivalent to $\pi_{mk}$ \cite{bs,bebo}.

\subsection{Infinitesimal representations}

\subsubsection{Holomorphic realizations}

Although it is a non-reductive algebraic group, $G^{J}$ can be considered as a
group of Harish-Chandra type \cite{bs}. Then the complexification of $\hat
{\pi}(\mathfrak{g}^{J})$ is the direct sum of vector spaces $\hat{\pi
}(\mathfrak{g}_{\mathbb{C}}^{J})={\mathfrak{p}}_{+}+{\mathfrak{t+p}}_{-\text{
}}$with ${\mathfrak{p}}_{-}={\mathfrak{p}}_{+}^{\dagger}$, $\left[
{\mathfrak{t,p}}_{\pm}\right]  \subset{\mathfrak{p}}_{\pm}$ , where
${\mathfrak{p}}_{+}=\left\langle a^{\dagger},K_{+}\right\rangle _{\mathbb{C}}%
$, and ${\mathfrak{t}}=\left\langle I\mathbf{,}K_{0}\right\rangle
_{\mathbb{C}}$. Consider the universal enveloping algebra $\mathcal{U(}%
{\mathfrak{p}}_{+})$ and the linearly independent elements $X_{1},\ldots
,X_{n}\in\mathcal{U(}{\mathfrak{p}}_{+})$. We now introduce a holomorphic
family of elements $E_{\mathbf{z}}=\exp(z_{1}X_{1}+\cdots+z_{n}X_{n})$, where
$\mathbf{z}=(z_{1},\ldots,z_{n})\in\Omega$ and $\Omega$ is an open subset of
$\mathbb{C}^{n}$. Suppose that there exists a cyclic vector $\Phi_{0}%
\in\mathcal{D}$ with ${\mathfrak{p}}_{-}\Phi_{0}=\{0\}$, ${\mathfrak{t}}%
\Phi_{0}=\left\langle \Phi_{0}\right\rangle _{\mathbb{C}}$ and $\left\Vert
\Phi_{0}\right\Vert =1$. This assumption is supported by \cite{bs}. Consider
now the holomorphic vectors $\Phi_{\mathbf{z}}=E_{\mathbf{z}}\Phi_{0}$, where
$\mathbf{z}\in\Omega$.

The map $T:\mathcal{H}\rightarrow\mathcal{H}_{\mathrm{hol}}(\Omega\mathcal{)}$
is defined by $T(\varphi)\left(  \mathbf{z}\right)  =\left\langle
\Phi_{\mathbf{\bar{z}}},\varphi\right\rangle $, $\varphi\in\mathcal{H}$,
$\mathbf{z\in}\Omega$, where $\mathcal{H}_{\mathrm{hol}}(\Omega\mathcal{)}$ is
the space of all $\mathbb{C}$-valued holomorphic functions on $\Omega$.
Consider an inner product on $T(\mathcal{H)}$ such that $T$ is unitary. We
will to obtain an infinitesimal irreducible unitary representation $T$%
$\hat{\pi}T^{-1}$on $T(\mathcal{H)}$ and the explicit form of the
corresponding basis differential operators. Let
\begin{equation}
\Phi_{\alpha w\tau}=\exp(\alpha a^{\dagger}+\frac{w}{2}a^{\dagger2}+\tau
W_{+})\Phi_{0},\label{holo}%
\end{equation}
and $\Psi_{\alpha w\tau}=\left\Vert \Phi_{\alpha w\tau}\right\Vert ^{-1}%
\Phi_{\alpha w\tau}$. Then the Bargmann coherent states \cite{per}, the
$SL(2,\mathbb{R)}$ coherent states \cite{per}, and the Perelomov coherent
states for $G^{J}$ \cite{jac1} are realized by the holomorphic vectors
$\Psi_{\alpha00}$ with $\alpha\in\mathbb{C}$, $\Psi_{00\tau}$ with $\tau
\in\mathbb{C}$, and $\Psi_{\alpha ww}$ with $\left(  \alpha,w\right)
\in\mathbb{C\times}\mathbb{D}$, respectively.

\subsubsection{ Schr\"{o}dinger-Weil representation of $G^{J}$}

Let $\mathbb{\mathcal{H}}=L^{2}(\mathbb{R)=\mathcal{H}}_{0}$ with the Schwartz
space $\mathcal{D=S}(\mathbb{R)}$ \cite{bs}. Let $\pi=\pi_{\mathrm{SW}}^{m}$
be the Schr\"{o}dinger-Weil representation of $G^{J}$. The basis differential
operators of $\hat{\pi}_{\mathrm{SW}}^{m}(\mathfrak{g}^{J}{\mathfrak{)}}$ can
be written as \cite{bs}:%
\begin{equation}
P=\frac{\mathrm{d}}{\mathrm{d}q},\ Q=2\mathrm{i}\mu q,\ R=\mathrm{i}\mu
I,\ F=\mathrm{i}\mu q^{2},\ G=\frac{\mathrm{i}}{4\mu}\frac{\mathrm{d}^{2}%
}{\mathrm{d}q^{2}},\ H=q\frac{\mathrm{d}}{\mathrm{d}q}+\frac{1}{2}%
I.\label{swdiff}%
\end{equation}
Then $2K_{-}=a^{2}$. In quantum mechanics, $\hbar=(2\mu^{2})^{-1}$ is the
Plank constant, $q$ is the position operator, $p\mathbf{=-}\mathrm{i}%
\hbar\mathrm{d}/\mathrm{d}q$ is the momentum operator, and $a$ and
$a^{\dagger}$ are the Fock annihilation and creation operators, respectively
\cite{per}. The vacuum vector $\varphi_{0}=\left(  2\left\vert \mu\right\vert
/\pi\right)  ^{1/4}\exp\left(  -\left\vert \mu\right\vert x^{2}\right)  $ and
the number vectors$\ \varphi_{n}=(n!)^{-1/2}(a^{\dagger})^{n}\varphi_{0}$%
,$\ $where $n\in\mathbb{N}$, form a complete orthonormal basis of analytic
vectors for the Hilbert space $\mathbb{\mathcal{H}}_{0}$. \noindent Consider
the map $T_{\mathrm{B}}(\varphi):\mathcal{H}_{0}\rightarrow\mathcal{H}%
_{\mathrm{hol}}(\mathbb{C}\mathcal{)}$ defined by $T_{\mathrm{B}}%
(\varphi)\left(  \alpha\right)  =\left\langle \Phi_{\bar{\alpha}00}%
,\varphi\right\rangle $, $\varphi\in\mathcal{H}_{0}$, $\alpha\in\mathbb{C}$,
where $\Phi_{0}=\varphi_{0}$. The polynomials $f_{\mathrm{B}n}^{\ }%
=T_{\mathrm{B}}(\varphi_{n})$, defined by $f_{\mathrm{B}n}^{\ }(z)=(n!)^{-1/2}%
z^{n}$, where $n\in\mathbb{N}$, form \ a complete orthonormal basis of
analytic vectors in the Hilbert space $T_{\mathrm{B}}(\mathcal{H}%
_{0})=L_{\mathrm{hol}}^{2}(\mathbb{C},\mu_{\mathrm{B}})$ with the Bargmann
measure given by $\mathrm{d}\mu_{\mathrm{B}}=\pi^{-1}\exp\left(  -\left\vert
z\right\vert ^{2}\right)  \mathrm{d}^{2}z$. The geometric quantization of
$\pi_{\mathrm{SW}}^{m}$ and $\pi_{\mathrm{B}}=T_{\mathrm{B}}\pi_{\mathrm{SW}%
}^{m}T_{\mathrm{B}}^{-1}$ is presented in \cite{gsa}. \ \noindent Consider now
the map $T_{0}:\mathcal{H}_{0}\rightarrow\mathcal{H}_{\mathrm{hol}%
}(\mathbb{C\times}\mathbb{D}\mathcal{)}$ defined by $T_{0}(\varphi)\left(
\alpha,w\right)  =\left\langle \Phi_{\bar{\alpha}\bar{w}0},\varphi
\right\rangle $, $\varphi\in\mathcal{H}_{0}$, $\left(  \alpha,w\right)
\in\mathbb{C\times}\mathbb{D}$, where $\Phi_{0}=\varphi_{0}$. The polynomials
$f_{n}^{\ }=T_{0}(\varphi_{n})$, where $n\in\mathbb{N}$, form \ a complete
orthonormal basis of analytic vectors for the Hilbert space $T_{0}%
(\mathcal{H}_{0})$. We have

\textbf{Proposition2.} \textit{a) The generating function of the polynomials
}$f_{n}$\textit{\ can be written as}
\begin{equation}
T_{\mathrm{B}}(\Phi_{\alpha w0})(z)=\exp\left(  \alpha z+\frac{1}{2}%
wz^{2}\right)  =\sum_{n\geq0}\frac{z^{n}}{\sqrt{n!}}f_{n}^{\ }(\alpha
,w).\ \label{generating}%
\end{equation}
b) \textit{Any solution} $f$ $\in\mathcal{H}_{\mathrm{hol}}(\mathbb{C\times
}\mathbb{D}\mathcal{)}$ \textit{of the equation} $\left(  \partial
^{2}/\partial\alpha^{2}-2\partial/\partial w\right)  f=0$ \textit{can be
written as }$f=\sum_{n\geq0}c_{n}f_{n}$, where $c_{n}\in\mathbb{C}$.
$\ $\textit{The map }$T_{0}$ \textit{is} \textit{non-surjective}.

\noindent c) \textit{Let} $\pi_{0}=T_{0}\pi_{\mathrm{SW}}^{m}T_{0}^{-1}$.
\textit{The basis differential operators can be expressed as}%
\begin{equation}
~\hat{\pi}_{0}(a)=\frac{\partial}{\partial\alpha},\ \hat{\pi}_{0}(K_{-}%
)=\frac{\partial}{\partial w},\ \hat{\pi}_{0}(K_{0})=\frac{1}{4}+\frac{\alpha
}{2}\frac{\partial}{\partial\alpha}+w\frac{\partial}{\partial w},\ \label{coh}%
\end{equation}%
\[
\hat{\pi}_{0}(a^{\dagger})=\alpha+w\frac{\partial}{\partial\alpha},\ \hat{\pi
}_{0}(K_{+})=\frac{1}{2}\alpha^{2}+\frac{w}{2}+\alpha w\frac{\partial
}{\partial w}+w^{2}\frac{\partial}{\partial w}.
\]

The Schr\"{o}dinger-Weil representation $\pi_{\mathrm{SW}}^{m}$, the Bargmann
representation $\pi_{\mathrm{B}}$, and the coherent state representation
$\pi_{\mathrm{0}}$ of $G^{J}$are unitarily equivalent.\ 

\subsubsection{Discrete series representations}

Let $k>1/2$.\emph{\ } Using Waldspurger's representation theory of
$Mp(2,\mathbb{R)}$ \cite{bs,wald} and \cite{sally}, we consider the
infinitesimal representations $\hat{\pi}_{k}$ of the positive discrete series
on the Hilbert space $\mathcal{H}_{k\ }$characterized by the normalized cyclic
vector $\phi_{0}$ and the complete orthonormal basis $\phi_{n^{\prime}n}%
^{[k]}$, where $n$, $n^{\prime}\in\mathbb{N}$, with $a\phi_{0}=0$, $~K_{-}%
\phi_{0}=0$, $~2K_{0}\phi_{0}=k\phi_{0}$,$\ $and $\phi_{n^{\prime}n}%
^{[k]}=C_{n^{\prime}n}^{[k]}(a^{\dagger})^{n^{\prime}}\left(  {W}\mathbf{_{+}%
}\right)  ^{n}\phi_{0}$. Here $C_{n^{\prime}n}^{[k]}=\left[  n!n^{\prime
}!(k-1/2)_{n}\right]  ^{-1/2}$ with $(x)_{n}=\Gamma(n+x)/\Gamma(x)$. Consider
now the map $T_{k}(\varphi):\mathcal{H}_{k}\rightarrow\mathcal{H}%
_{\mathrm{hol}}(\mathbb{C}\mathcal{)}$ defined by $T_{k}(\varphi)\left(
z,\zeta\right)  =\left\langle \Phi_{\bar{z}0\bar{\zeta}},\varphi\right\rangle
$, $\varphi\in\mathcal{H}_{k}$, where $\Phi_{0}=\phi_{0}$ and $z$, $\zeta
\in\mathbb{C}$ with $|\zeta|<1$. \noindent\ 

The polynomials $f_{n^{\prime}n}^{\ [k]}=$ $T_{k}(\phi_{n^{\prime}n}^{[k]})$,
where $n$, $n^{\prime}\in\mathbb{N}$, form \ a complete orthonormal basis of
analytic vectors for the Hilbert space $T_{k}(\mathcal{H}_{k})$. Here
$f_{n^{\prime}n}^{\ [k]}(z,\zeta)=D_{n^{\prime}n}^{[k]}z^{n^{\prime}}\zeta
^{n}$, where $D_{n^{\prime}n}^{[k]}=(n!n^{\prime}!)^{-1/2}(k-1/2)_{n\ }%
^{\ 1/2}$. Consider $\hat{\sigma}_{k}=T_{k}\hat{\pi}_{k}T_{k}^{-1}$. The basis
differential operators can be written as
\begin{align}
\hat{\sigma}_{k}(a)  &  =\frac{\partial}{\partial z},\ \hat{\sigma}%
_{k}(a^{\dagger})=z,\ \hat{\sigma}_{k}(K_{0})=\frac{k}{2}+\frac{1}{2}%
z\frac{\partial}{\partial z}+\zeta\frac{\partial}{\partial\zeta}%
,\label{discretediff}\\
\hat{\sigma}_{k}(K_{-})  &  =\frac{1}{2}\frac{\partial^{2}}{\partial z^{2}%
}+\frac{\partial}{\partial\zeta},\ \hat{\sigma}_{k}(K_{+})=\frac{z^{2}}%
{2}+(k-\frac{1}{2})\zeta+\zeta^{2}\frac{\partial}{\partial\zeta}.\ \nonumber
\end{align}
Let $\mathrm{d}\mu_{k}=(k-3/2)\pi^{-2}\exp\left(  -\left\vert z\right\vert
^{2}\right)  (1-\left\vert \zeta\right\vert ^{2})^{k-5/2}\mathrm{d}%
^{2}z\ \mathrm{d}^{2}\zeta$. If $k>3/2$, then $T_{k}(\mathcal{H}%
_{k})=L_{\mathrm{hol}}^{2}(\mathbb{C\times D},\mu_{k})$.There is an analytic
continuation of\emph{\ }$\hat{\pi}_{k}$ in the limit $k\rightarrow3/2$ .

\section{Squeezed states for the Jacobi group}

\subsection{Matrix elements}

We consider the \textit{squeezed operator }$T(\alpha,w)=D\left(
\alpha\right)  S\left(  w\right)  $ for\textit{\ }$G^{J}$, where the
\textit{displacement operator }$D\left(  \alpha\right)  $ and the \textit{the
squeezed operator }$S\left(  w\right)  $ are unitary operators defined by
\cite{per}:%
\begin{equation}
D(\alpha)=\exp(\alpha a^{\dagger}-\bar{\alpha}a)=\exp(-\frac{1}{2}|\alpha
|^{2})\exp(\alpha a^{\dagger})\exp(-\bar{\alpha}a),\ \label{displacement}%
\end{equation}%
\begin{equation}
S\left(  w\right)  =\exp(wK_{+})\exp(\eta K_{0})\exp(-\bar{w}K_{-}%
),\ \label{squeezed}%
\end{equation}
where $\alpha\in\mathbb{C}$, $w\in\mathbb{C}$,$\ |w|<1$, and $\eta
=\ln(1-|w|^{2})$. The matrix elements of $D(\alpha)$ are given by Schwinger's
formula in terms of Laguerre polynomials \cite{per}. The matrix elements of
$S\left(  w\right)  $ for the Schr\"{o}dinger-Weil representation can be
written as associated Legendre functions \cite{per}. Moreover, the matrix
elements of $S\left(  w\right)  $ for discrete series representations of
$G^{J}$can be expressed in terms of hypergeometric polynomials \cite{sally}:
\begin{equation}
\left\langle \phi_{0n^{\prime}}^{[k]}\right\vert S\left(  w\right)  \left\vert
\phi_{0n}^{[k]}\right\rangle \!=\!\frac{\lambda_{kn}w^{s}}{\lambda
_{kn^{\prime}}s!}\!\!\left(  \!1\!-|w|^{2}\right)  ^{h}\!F\left(
\!-n,\!s+n+2h;s+1\!;\!|w|^{2}\right)  ,\label{squeezedmatrix}%
\end{equation}
where $n$, $n^{\prime}\in\mathbb{N}$, $k>1/2$,\ $h=(2k-1)/4$, $\lambda
_{kc}=\left[  c!\Gamma(2h+c)/\Gamma(2h)\right]  ^{-1/2}$ for $c=n,n^{\prime}$,
and $s=n^{\prime}-n\geq0.$

We now introduce the \textit{squeezed state vectors }$T(\alpha,w)\varphi$
\ for the Jacobi group $G^{J}$, where $\varphi\in\mathcal{D}$ and $\left\Vert
\varphi\right\Vert =1$. The standard coherent states, squeezed states,
displaced number states, squeezed number states, and displaced squeezed number
states \cite{ni} are realized by $T(\alpha,0)\varphi_{0}$, $T(0,w)\varphi_{0}%
$, $T(\alpha,0)\varphi_{n}$, $T(0,w)\varphi_{n}$, and $T(\alpha,w)\varphi_{n}%
$, respectively \cite{ni}. The squeezed states for\textit{\ } $G^{J}$ can be
in particular the coherent states introduced by Schr\"{o}dinger \cite{sch},
the squeezed states considered by Kennard \cite{ken}, and the displaced
squeezed number states of Husimi \cite{husi}.

We now introduce the notation $\hat{A}=S\left(  -w\right)  D\left(
-\alpha\right)  AD\left(  \alpha\right)  S\left(  w\right)  $ and
$r=(1-|w|^{2})^{-1/2}$. We obtain $\hat{a}=(\widehat{a^{\dagger}})^{\dagger
}=r\left(  a+wa^{\dagger}\right)  +\alpha I$ and
\[
\hat{K}_{-}=\hat{K}_{+}^{\dagger}=r^{2}\left(  K_{-}+2wK_{0}+w^{2}%
K_{+}\right)  +r\alpha\left(  a+wa^{\dagger}\right)  +\frac{\alpha^{2}}%
{2}I\mathbf{,}%
\]%
\[
\hat{K}_{0}=r^{2}\left[  \bar{w}K_{-}+\left(  1+|w|^{2}\right)  K_{0}%
+wK_{+}\right]  +r\Re\alpha\left(  a^{\dagger}+\bar{w}a\right)  +\frac
{|\alpha|^{2}}{2}I.
\]
If $A=a^{\dagger n}K_{+}^{m}K_{0}^{s}K_{-}^{m^{\prime}}a^{n^{\prime}}$, then
$\hat{A}=\hat{a}^{\dagger n}\hat{K}_{+}^{m}\hat{K}_{0}^{s}\hat{K}%
_{-}^{m^{\prime}}\hat{a}^{n^{\prime}}$. Using the preceding results, we can
obtain the expectation values of any polynomial operator in infinitesimal
generators of $G^{J}$.

\subsection{Uncertainty relations}

Consider the Schr\"{o}dinger inequality \cite{sch2} $\sigma_{AA}\sigma
_{BB}\geq\sigma_{AB}^{2}+\left\vert \left\langle [A,B]\right\rangle
\right\vert ^{2}/4$ for the selfadjoint operators $A$ and $B$ , where
$\sigma_{AB}=\left\langle AB+BA\right\rangle /2-\left\langle A\right\rangle
\left\langle B\right\rangle $ and $\sigma_{CC}=\left\langle C^{2}\right\rangle
-\left\langle C\right\rangle ^{2}$ for $C=A,B$. Here $\langle\,\rangle$ means
the expectation value with respect to the state vector $\Phi\in\mathcal{D}$.
Consider now $\Phi=T(\alpha,w)\varphi_{n}$ or $\Phi=T(\alpha,w)\phi
_{nn^{\prime}}^{[k]}$. Let $u_{\pm}=r^{2}\left(  1\pm w\right)  \left(
1\pm\bar{w}\right)  $ and $n_{0}=n+1/2$. We obtain $\sigma_{qq}=n_{0}\hbar
u_{+}$,$\ ~\sigma_{pp}=n_{0}\hbar u_{-}$, and $\sigma_{pq}=2n_{0}\hbar
r^{2}\Im w$.

We have $\sigma_{qq}\sigma_{pp}=\sigma_{pq}^{2}+\hbar^{2}/4$ for the squeezed
states $T(\alpha,w)\varphi_{0}$ and $T(\alpha,w)\phi_{0}$. Moreover,
$\sigma_{qq}\sigma_{pp}=\hbar^{2}/4$ for the coherent states $T(\alpha
,0)\varphi_{0}$ and $T(\alpha,0)\phi_{0}$. \ Evidently, $\sqrt{\sigma
_{qq}\sigma_{pp}}\geq n_{0}\hbar$, and we have \textit{squeezing} in the
region $2n_{0}u_{+}<1$, described by the open disk \ $\left\vert \left(
2n_{0}+1\right)  w+2n_{0}\right\vert \,<1$, $w\in\mathbb{C}$.

\subsection{Mandel's parameter}

Consider the Mandel parameter $Q=\left\langle (\Delta N)^{2}\right\rangle
/\left\langle N\right\rangle -1$, where $N=a^{\dagger}a$ is the number
operator. We obtain
\begin{equation}
Q\left(  \alpha.w\right)  =\frac{(4n_{0}{}^{2}+3)|w|^{2}+4n_{0}|\alpha\bar
{w}+\bar{\alpha}|^{2}(1-|w|^{2})}{2n_{0}(1-|w|^{4})+\left(  2|\alpha
|^{2}-1\right)  (1-|w|^{2})^{2}}-1. \label{mandel}%
\end{equation}
We have $Q\left(  \alpha,0\right)  =n\left(  2|\alpha|^{2}-1\right)  \left(
n+|\alpha|^{2}\right)  ^{-1}$. Then $Q\left(  \alpha,0\right)  =0$ for $n=0$
or $|\alpha|=1/\sqrt{2\text{ }}$. Moreover, $Q\left(  0,w\right)  =0$ for
\begin{equation}
|w|^{2}=\frac{1}{2\left(  2n_{0}+1\right)  }\left(  \sqrt{16n_{0}^{4}%
+24n_{0}^{2}-3}-4n_{0}^{2}-1\right)  . \label{mandel2}%
\end{equation}
The preceding formulas are compatible with \cite{oli,quant,ni}.

\end{document}